\documentclass[12pt]{article}%
\usepackage{amsmath}
\usepackage{amscd}
\usepackage{amsfonts}
\usepackage{amssymb}%
\newtheorem{theorem}{Theorem}[section]
\newtheorem{lemma}[theorem]{Lemma}

\newtheorem{definition}[theorem]{Definition}
\newtheorem{corollary}[theorem]{Corollary}
\newtheorem{remark}[theorem]{Remark}

\newenvironment{proof}{\bf Proof. \rm}{$\Box$}
\newcommand{\be}{\begin{equation}}
\newcommand{\ee}{\end{equation}}

\begin{document}

\title{Extensions and Dilations for $C^{*}$-dynamical Systems}
\author{Paul S. Muhly\thanks{Supported in part by grants from the National Science
Foundation (DMS-0355443) and from the U.S.-Israel Binational Science
Foundation.}\\Department of Mathematics\\University of Iowa\\Iowa City, IA 52242\\e-mail: pmuhly@math.uiowa.edu
\and Baruch Solel\thanks{Supported in part by the U.S.-Israel Binational Science
Foundation and by the Fund for the Promotion of Research at the Technion.}\\Department of Mathematics\\Technion\\32000 Haifa, Israel\\e-mail: mabaruch@techunix.technion.ac.il}
\date{}
\maketitle

\section{Introduction}

Throughout this note, $A$ will denote a $C^{\ast}$-algebra with unit
$\mathbf{1}$ and $\alpha$ will denote an injective, unital endomorphism of $A$.

\begin{definition}
A \emph{(contractive) covariant representation} of the pair $(A,\alpha)$ is a
pair $(\pi,T)$ consisting of a $C^{\ast}$-representation $\pi$ of $A$ on a
Hilbert space $H$ and a contraction operator $T$ in $B(H)$ such that%
\begin{equation}
T\pi(\alpha(a))=\pi(a)T\text{,}\label{Cov}%
\end{equation}
for all $a\in A$. If $T$ is an isometry, then we say that $(\pi,T)$ is
\emph{isometric}, while if $T$ is a coisometry, we say $(\pi,T)$ is
\emph{coisometric.}
\end{definition}

Our primary objective here is to prove the following theorem and corollary.

\begin{theorem}
\label{coisoext} If $(\pi,T)$ is a contractive covariant representation of
$(A,\alpha)$ on a Hilbert space $H$, then there exists a coisometric covariant
representation $(\rho,V)$ on a Hilbert space $K$ containing $H$ that extends
$(\pi,T)$. That is, $\rho(a)H\subseteq H$ and $\rho(a)|H=\pi(a)$ for all $a\in
A$, while $VH\subseteq H$ and $V|H=T$.
\end{theorem}

\begin{corollary}
\label{dilation}If $(\pi,T)$ is a contractive covariant representation of
$(A,\alpha)$ on a Hilbert space $H$, then there exists a covariant
representation $(\sigma,U)$ of $(A,\alpha)$ on a Hilbert space $K$ containing
$H$ such that $U$ is \emph{unitary} and such that $T^{n}=PU^{n}|H$ for all
$n\geq0$, where $P$ denotes the projection of $K$ onto $H$.
\end{corollary}

Theorem \ref{coisoext} is a special case of Theorem 5.10 of \cite{MS98} while
Corollary \ref{dilation} is a special case of Theorem 5.22 in \cite{MS98}. Our
exposition is designed to give self-contained, elementary proofs of these
results which avoid the technology employed in \cite{MS98}. That is, we avoid
the formal use of the theory of $C^{\ast}$-correspondences, which, in a sense,
are the central objects of \cite{MS98}. We hope that the exposition given here
will aid the interested reader in his or her efforts to understand the results
of \cite{MS98}.

In contrast to the situation for single contraction operators on Hilbert
space, neither the coisometric extension $(\rho,V)$ of $(\pi,T)$ in Theorem
\ref{coisoext} nor the dilation $(\sigma,U)$, with $U$ unitary, in Corollary
\ref{dilation} is unique in general. (By ``unique'', here, we really mean
``unique up to unitary equivalence''.) Indeed, unless $\alpha$ is an
automorphism, where it is possible to obtain uniqueness by imposing minimality
of $(\rho,V)$ and $(\sigma,U)$ as is done in \cite{MM83}, $(\rho,V)$ and
$(\sigma,U)$ cannot be chosen uniquely. An additional objective of this note,
then, is to explain the lack of uniqueness and to discuss a way to organize
the coisometric extensions and dilations of a contractive covariant
representation of $(A,\alpha)$ in a special, but important, situation.

To see why a dilation $(\sigma,U)$ of $(\pi, T)$, with $U$ unitary, might not
be unique unless $\alpha$ is an automorphism, observe that the covariance
equation, equation (\ref{Cov}), implies that
\[
\sigma(\alpha(a))=U^{\ast}\sigma(a)U\text{,}%
\]
for all $a\in A$. Thus, as Stacey develops in \cite{pS93}, the representation
$\sigma$ extends to the algebra $A_{\infty}$, which is the inductive limit of
the inductive system%
\[
A\overset{\alpha}{\rightarrow}A\overset{\alpha}{\rightarrow}A\overset{\alpha
}{\rightarrow}\cdots
\]
built from $A$ and $\alpha$. So, to construct $(\sigma,U)$ from $(\pi,T)$ we
are faced with the problem of extending $\pi$ to $A_{\infty}$ and this
extension problem does not have a unique solution. Of course, in the process
of extending $\pi$ to $A_{\infty}$ to obtain $\sigma$ we really have to
\textquotedblleft inter-leave\textquotedblright\ the construction of $\sigma$
with the construction of $U$. One might hope that this would cut down on the
possibilities for $\sigma$ leading to an essentially unique dilation
$(\sigma,U)$. However, it doesn't. As we will see, the same problem arises
when building $(\rho,V)$. Indeed, $(\rho,V)$ lies at the center of the
uniqueness problem.

We note, too, that our standing hypotheses, that $A$ is unital and that
$\alpha$ preserves the unit and is injective, are made, fundamentally, to
avoid complications with inductive limits. This may not be evident from the
proofs, which proceed at a very elementary level, but closer analysis reveals
that once these hypotheses are relaxed, then difficulties begin to arise. The
presentation in \cite{MS98} is designed to address these - and to extend to a
much broader context.

\section{Proofs}

The proof of Theorem \ref{coisoext} is based on two lemmas. The first is

\begin{lemma}
\label{HB} Given a representation $\pi:A\rightarrow B(H)$, there is a Hilbert
space $K$, an isometry $W:H\rightarrow K$ and a representation $\rho$ of $A$
in $B(K)$ such that
\begin{equation}
W^{\ast}\rho(\alpha(a))W=\pi(a),\label{Ext}%
\end{equation}
for all $a\in A.$
\end{lemma}

Note: Equation (\ref{Ext}) implies that $WW^{\ast}$ commutes with $\rho
(\alpha(A))$. This fact will play an important role throughout the
computations that use this lemma.

\begin{proof}
If $\pi$ has a unit cyclic vector $\xi$, say, let $\omega_{0}$ be the state on
$\alpha(A)$ defined by the formula, $\omega_{0}(\alpha(a))=(\pi(a)\xi,\xi)$.
Apply the Hahn-Banach theorem to extend $\omega_{0}$ to a state $\omega$ on
$A$ and let $(\rho,K)$ be the GNS representation of $A$ determined by $\omega
$. Then $W$ is defined by setting $W\pi(a)\xi=\rho(\alpha(a))\xi$. In general,
we may write $\pi$ as a direct sum of cyclic representations and apply this
argument to each summand.
\end{proof}

\bigskip

Of course the Hahn-Banach theorem introduces an arbitrariness to the extension
$\rho$ of $\pi$ and the operator $W$ that cannot be avoided in general; it is
the source of non-uniqueness in the theory. However, it can be controlled if
we have a \textquotedblleft transfer operator\textquotedblright\ for $\alpha$
at our disposal.

\begin{definition}
\label{Transfer}A \emph{transfer operator} for $\alpha$ is a completely
positive left inverse of $\alpha$; i.e., a completely positive map
$\tau:A\rightarrow A$ such that $\tau\circ\alpha(a)=a$ for all $a\in A$.
\end{definition}

The idea of introducing transfer operators into the study of endomorphisms of
$C^{*}$-algebras is due to R. Exel \cite{rE03}. On the face of it, our
definition is a bit different from his. However, he works in a more general
setting than ours and it is not hard to see that under our hypotheses that
$\alpha$ is injective and unital, his definition coincides with ours (see his
Proposition 2.6 in particular). As Exel explains, in general a transfer
operator need not exist for $\alpha$, and when one does exist, it need not be
unique. In fact, as he shows, transfer operators are just as plentiful as
conditional expectations of $A$ onto the range of $\alpha$. To see this, note
that if $\tau$ is a transfer operator for $\alpha$, then $E:=\alpha\circ\tau$
is a conditional expectation onto the range of $\alpha$. Indeed, $E$ evidently
is a completely positive, unital map with range $\alpha(A)$ and the
computation,
\[
E^{2}=(\alpha\circ\tau)\circ(\alpha\circ\tau)=\alpha\circ(\tau\circ
\alpha)\circ\tau=\alpha\circ\tau=E\text{,}%
\]
completes the proof. Conversely, if $E$ is a conditional expectation onto the
range of $\alpha$, then ${\alpha}^{-1}\circ E$ is a transfer operator for
$\alpha$. In the commutative setting, transfer operators are the subject of an
active area of study because of their importance in statistical mechanics,
Markov processes and the general theory of irreversible dynamical systems (see
\cite{vB00}).

\begin{remark}
\label{TransferUnique}If $\tau$ is a transfer operator for $\alpha$ and if
$\pi$ is a representation of $A$ on the Hilbert space $H$, then $\pi\circ\tau$
is a unital completely positive map from $A$ to $B(H)$. Further, if $\rho$ is
the minimal Stinespring dilation of $\pi\circ\tau$, mapping $A$ to $B(K)$, and
if $W:H\rightarrow K$ is the isometry that comes in Stinespring's theorem,
then for all $a\in A$,
\[
W^{\ast}\rho(\alpha(a))W=\pi\circ\tau(\alpha(a))=\pi(a).
\]
By the minimality assumption the pair $(\rho,W)$ is uniquely determined up to
unitary equivalence. (See the Remark after \cite[Theorem 1.1.1]{wA69}.)

If $\alpha$ has a transfer operator $\tau$ and if $\pi:A\rightarrow B(H)$ is a
representation, then we shall say that the triple $(\rho,W,K)$ obtained from
the minimal Stinespring dilation of the completely positive map $\pi\circ\tau
$, mapping $A$ to $B(H)$, is \emph{the} \emph{extension of }$\pi$
\emph{adapted to }$\tau$.
\end{remark}

\begin{lemma}
\label{2step}Given $(\pi,T)$ acting on $H$, choose $\rho$, $W$ and $K$ as in
Lemma \ref{HB}. Define the following objects: $\Delta_{\ast}:=(I-TT^{\ast
})^{1/2}$ (note that $\Delta_{\ast}$ commutes with $\pi(A)$); $\mathcal{D}%
_{\ast}:=\overline{\rho(A)W\Delta_{\ast}H}\subseteq K$; $D_{\ast}%
:=\Delta_{\ast}W^{\ast}|\mathcal{D}_{\ast}$; and the representation
$\hat{\pi}:A\rightarrow B(\mathcal{D}_{\ast})$, where $\hat{\pi}(a):=\rho
(a)|\mathcal{D}_{\ast}$. Then $\left(  \left[
\begin{array}
[c]{cc}%
\pi & \\
& \hat{\pi}
\end{array}
\right]  ,\left[
\begin{array}
[c]{cc}%
T & D_{\ast}\\
0 & 0
\end{array}
\right]  \right)  $, acting on $H\oplus\mathcal{D}_{\ast}$, gives a
contractive covariant representation of $A$ such that $\left[
\begin{array}
[c]{cc}%
T & D_{\ast}\\
0 & 0
\end{array}
\right]  $ is a partial isometry and which, when restricted to $H$, gives
$(\pi,T)$.
\end{lemma}

\begin{proof}
This is a simple matrix computation:%
\begin{multline*}
\left[
\begin{array}
[c]{cc}%
T & D_{\ast}\\
0 & 0
\end{array}
\right]  \left[
\begin{array}
[c]{cc}%
\pi(\alpha(a)) & \\
& \hat{\pi}(\alpha(a))
\end{array}
\right]  \\
=\left[
\begin{array}
[c]{cc}%
T\pi(\alpha(a)) & D_{\ast}\hat{\pi}(\alpha(a))\\
0 & 0
\end{array}
\right]  \\
=\left[
\begin{array}
[c]{cc}%
T\pi(\alpha(a)) & \Delta_{\ast}W^{\ast}\rho(\alpha(a))WW^{\ast}%
|\mathcal{D}_{\ast}\\
0 & 0
\end{array}
\right]  \\
=\left[
\begin{array}
[c]{cc}%
\pi(a)T & \Delta_{\ast}W^{\ast}\rho(\alpha(a))WW^{\ast}|\mathcal{D}_{\ast
}\\
0 & 0
\end{array}
\right]  \\
=\left[
\begin{array}
[c]{cc}%
\pi(a)T & \Delta_{\ast}\pi(a)W^{\ast}|\mathcal{D}_{\ast}\\
0 & 0
\end{array}
\right]  \\
=\left[
\begin{array}
[c]{cc}%
\pi(a)T & \pi(a)\Delta_{\ast}W^{\ast}|\mathcal{D}_{\ast}\\
0 & 0
\end{array}
\right]  \\
=\left[
\begin{array}
[c]{cc}%
\pi(a)T & \pi(a)D_{\ast}\\
0 & 0
\end{array}
\right]  \\
=\left[
\begin{array}
[c]{cc}%
\pi(a) & \\
& \hat{\pi}(a)
\end{array}
\right]  \left[
\begin{array}
[c]{cc}%
T & D_{\ast}\\
0 & 0
\end{array}
\right]
\end{multline*}

\end{proof}

\bigskip
As is customary in the dilation theory of single operators, $\Delta_*$ is called the \emph{defect operator} of $T^*$ and $\mathcal{D}_{*}$ is called the associated \emph{defect space}.
\bigskip

\textbf{Proof of Theorem \ref{coisoext}.} Lemma \ref{2step} is the
\textquotedblleft zeroth\textquotedblright\ step in an inductive construction.
We apply Lemma \ref{HB} again to $\hat{\pi}$, which is a representation of
$A$ on $\mathcal{D}_{\ast}$, to obtain a Hilbert space $K_{1}$, an isometric
embedding $W_{1}:\mathcal{D}_{\ast}\overset{}{\longrightarrow}K_{1}$ and a
rep. $\pi_{1}:A\rightarrow B(K_{1})$ such that $W_{1}^{\ast}\pi_{1}%
(\alpha(a))W_{1}=\hat{\pi}(a)$, for all $a\in A$. Observe that $W_{1}%
W_{1}^{\ast}$ commutes with $\pi_{1}(\alpha(A))$. Set

\begin{enumerate}
\item $\mathcal{D}_{1\ast}:=\overline{\pi_{1}(A)W_{1}\mathcal{D}_{\ast}%
}\subseteq K_{1}$,

\item $D_{1\ast}:=W_{1}^{\ast}$, and

\item $\hat{\pi}_{1}:A\rightarrow B(\mathcal{D}_{1\ast})$, $\hat{\pi}_{1}%
(a):=\pi_{1}(a)|\mathcal{D}_{1\ast}$
\end{enumerate}

Note that $W_{1}\mathcal{D}_{\ast}\subseteq\mathcal{D}_{1\ast}$ and that while
$\hat{\pi}_{1}(\alpha(A))W_{1}\mathcal{D}_{\ast}\subseteq W_{1}\mathcal{D}_{\ast
}$, $\hat{\pi}_{1}(A)\mathcal{D}_{\ast}$ need not be contained in $\mathcal{D}%
_{\ast}$.

Inductively, we obtain sequences $\{\mathcal{D}_{k\ast}\}_{k\geq1}$,
$\{W_{k}\}_{k\geq1}$, $\{\pi_{k}\}_{k\geq1}$, $\{\hat{\pi}_{k}\}_{k\geq1}$, \ and
$\{D_{k}\}_{k\geq1}$, where for $k>1$, $\mathcal{D}_{k\ast}$ and $K_{k}$ are
Hilbert spaces, $W_{k}:\mathcal{D}_{k-1\ast}\overset{}{\longrightarrow}$
$K_{k}$ is an isometry, $\pi_{k}:A\rightarrow B(K_{k})$ is a $C^{\ast}%
$-representation, and $\hat{\pi}_{k}:A\rightarrow B(\mathcal{D}_{k\ast})$ is a
$C^{\ast}$-representation such that the equations%
\[
W_{k}^{\ast}\pi_{k}(\alpha(a))W_{k}=\hat{\pi}_{k-1}(a)\text{,}%
\]%
\[
\mathcal{D}_{k\ast}:=\overline{\pi_{k}(A)W_{k}\mathcal{D}_{(k-1)^\ast}}\subseteq
K_{k-1}\text{,}%
\]%
\[
D_{k\ast}:=W_{k}^{\ast}\text{,}%
\]
and%
\[
\hat{\pi}_{k}:A\rightarrow B(\mathcal{D}_{k\ast}),\ \hat{\pi}_{k}(a):=\pi
_{k}(a)|\mathcal{D}_{k\ast}\text{,}%
\]
are satisfied.

On $H\oplus\mathcal{D}_{\ast}\oplus\mathcal{D}_{1\ast}\oplus\mathcal{D}%
_{2\ast}\oplus\cdots$ set%
\[
\rho:=\left[
\begin{array}
[c]{cccccc}%
\pi &  &  &  &  & \\
& \hat{\pi} &  &  &  & \\
&  & \hat{\pi}_{1} &  &  & \\
&  &  & \hat{\pi}_{2} &  & \\
&  &  &  & \ddots & \\
&  &  &  &  & \ddots
\end{array}
\right]
\]
and%
\[
V:=\left[
\begin{array}
[c]{cccccc}%
T & D_{\ast} &  &  &  & \\
& 0 & D_{1\ast} &  &  & \\
&  & 0 & D_{2\ast} &  & \\
&  &  & 0 & \ddots & \\
&  &  &  & \ddots & \ddots\\
&  &  &  &  & \ddots
\end{array}
\right]
\]
Then a straightforward calculation reveals that $(\rho,V)$ is a covariant
representation that extends $(\pi,T)$, where $V$ is a coisometry. $\Box$

\bigskip

\begin{remark}
\label{TransferUnique2}Of course, the multiple uses of Lemma \ref{HB}
contribute to the nonuniqueness of $(\rho,V)$. However, if there is a transfer
operator $\tau$ for $\alpha$ that is fixed in advance and if all the
sequences, $\{\mathcal{D}_{k\ast}\}_{k\geq1}$, $\{W_{k}\}_{k\geq1}$,
$\{\pi_{k}\}_{k\geq1}$, $\{\hat{\pi}_{k}\}_{k\geq1}$, \ and $\{D_{k}\}_{k\geq1}$,
are adapted to $\tau$ as in Remark \ref{TransferUnique}, then it is easy to
see that $(\rho,V)$ is unique up to unitary equivalence. If these sequences
are adapted to $\tau$, then we shall say that $(\rho,V)$\emph{ is adapted to
}$\tau$\emph{.}
\end{remark}

To prove Corollary \ref{dilation} we could dilate $(\rho,V)$ using
\cite[Theorem 3.3]{MS98}. As is proved there, in contrast to coisometric
extensions, the \emph{isometric dilation} of a contractive covariant representation is uniquely
determined (provided, of course, it is minimal in a well-known sense that we recapitulate in Theorem \ref{isodilate}). Further, as is shown in \cite[Theorem 2.18]{MS02}, the dilation
for $(\rho,V)$ will be coisometric and isometric. For completeness, we give a
proof of these results in our special situation which avoids the overhead of
the theory of $C^{\ast}$-correspondences.

\begin{theorem}
\label{isodilate}Let $(\pi,T)$ be a contractive covariant representation of
$(A,\alpha)$ on a Hilbert space $H$. Then there is a Hilbert space $K$
containing $H$ and an isometric covariant representation $(\eta,W)$ of
$(A,\alpha)$ on $K$ such that $\eta(a)H\subseteq H$ and $\eta(a)|H=\pi(a)$,
for all $a\in A$, and such that for all $n\geq0$, $T^{n}=PW^{n}|H$, where $P$
denotes the projection of $K$ onto $H$. Further:

\begin{enumerate}
\item $(\eta,W)$ is uniquely determined by $(\pi,T)$ up to unitary
equivalence if it is assumed (as may be arranged) that the smallest subspace
$K$ containing $H$ that is invariant under $W$ is $K$, i.e., if it is assumed that $(\eta,W)$ is \emph{minimal}; and

\item if $T$ is a coisometry, then so is $W$.
\end{enumerate}
\end{theorem}

\begin{proof}
The proof is to build the lower right-hand corner of the Schaeffer matrix for
the minimal unitary dilation of $T$ and to check that the representation $\pi$
can be extended to the Hilbert space of this dilation. For this purpose, let
$\Delta$ be the square root of $I-T^{\ast}T$, i.e., the defect operator of $T$. Thanks to the covariance
equation, equation (\ref{Cov}), $\Delta$ commutes with $\pi \circ \alpha$ and so the
closure of the range of $\Delta$, $\mathcal{D}$ - the defect space of $T$, reduces $\pi \circ \alpha$. If we let
$K:=H\oplus\mathcal{D}\oplus\mathcal{D}\oplus\cdots$ and on $K$ define $\eta$
and $W$ by the infinite matrices%
\[
\eta(a):=\left[
\begin{array}
[c]{ccccc}%
\pi &  &  &  & \\
& \pi\circ\alpha|\mathcal{D} &  &  & \\
&  & \pi\circ\alpha^{2}|\mathcal{D} &  & \\
&  &  & \pi\circ\alpha^{3}|\mathcal{D} & \\
&  &  &  & \ddots
\end{array}
\right]  \text{,}%
\]
$a\in A$, and%
\[
W:=\left[
\begin{array}
[c]{ccccc}%
T &  &  &  & \\
\Delta &  &  &  & \\
& I_{\mathcal{D}} &  &  & \\
&  & I_{\mathcal{D}} &  & \\
&  &  & \ddots &
\end{array}
\right]  \text{,}%
\]
then a straightforward calculation shows that $(\eta,W)$ is an isometric
covariant representation that dilates $(\pi,T)$. Further, it is evident that
$(\eta,W)$ is minimal in the sense that $K$ is the smallest subspace
containing $H$ that reduces $W$. It is also evident that if $T$ is a
coisometry, then so is $W$, i.e., in this event, $W$ is unitary. The
uniqueness assertion, 1., is immediate from the uniqueness of the minimal
isometric dilation of a contraction. See, for example, the proof on page 37 of
\cite{SzNF70}.
\end{proof}

\bigskip

\begin{definition}
The isometric dilation $(\eta,W)$ of a contractive covariant representation
$(\pi,T)$ constructed in Theorem \ref{isodilate} is called \emph{the minimal
isometric dilation of }$(\pi,T)$.\bigskip
\end{definition}

\textbf{Proof of Corollary \ref{dilation}. }As we indicated above, we simply
apply Theorem \ref{isodilate} to a coisometric extension $(\rho,V)$ of
$(\pi,T)$. The resulting covariant representation $(\sigma,U)$ will be
isometric and coisometric by the assertion 2. of the theorem. $\Box$ \bigskip

It may be helpful to have a matricial picture for $(\sigma,U)$. We follow the
notation developed above in the proof of Theorem \ref{coisoext}. Let $q_{0}$
be the projection of $\mathcal{D}_{\ast}$ onto $\mathcal{D}_{\ast}%
\ominus\overline{W\Delta_{\ast}H}$ and for $k\geq1$, let $q_{k}$ be the
projection of $\mathcal{D}_{k\ast}$ onto $\mathcal{D}_{k\ast}\ominus
W_{k}\mathcal{D}_{(k-1)\ast}$. Also, let $\Delta:=(I-T^{\ast}T)^{1/2}$ and
define
\begin{align*}
\mathcal{D}  &  :=\overline{\Delta H}\oplus q_{0}(\mathcal{D}_{\ast})\oplus
q_{1}(\mathcal{D}_{1\ast})\oplus\cdots\\
&  \subseteq H\oplus\mathcal{D}_{\ast}\oplus\mathcal{D}_{1\ast}\oplus
\mathcal{D}_{2\ast}\oplus\cdots
\end{align*}
- the Hilbert space for $(\rho,V)$. Here $\mathcal{D}$ is the defect space for $V$.  It contains the defect space of $T$, $\overline{\Delta H}$, as a summand. On $\mathcal{D}$ set
\[
\rho_{1}(a):=diag(\pi\circ\alpha(a),\hat{\pi}\circ\alpha(a),\hat{\pi}
_{1}\circ\alpha(a),\cdots).
\]
Since $\hat{\pi}_{k}(\alpha(A))W_{k}\mathcal{D}_{(k-1)\ast}\subseteq
W_{k}\mathcal{D}_{(k-1)\ast},$ $\rho_{1}$ is well defined. Observe that $\rho_1$ really is the restriction of $\rho \circ \alpha$ to $\mathcal{D}$. Let
\[X:=%
(-T^{\ast}W^{\ast}|\overline{(W\Delta_{\ast}H)})\oplus q_{0}:\mathcal{D}_{\ast
}\rightarrow\mathcal{D}.
\] On $\cdots\mathcal{D}_{2\ast}\oplus\mathcal{D}%
_{1\ast}\oplus\mathcal{D}_{\ast}\oplus H\oplus\mathcal{D}\oplus\mathcal{D}%
\oplus\cdots$, then, $U$ is represented matricially as{}%
\[
U=\left[
\begin{array}
[c]{ccccccccc}%
\ddots &  &  &  &  &  &  &  & \\
\ddots & 0 &  &  &  &  &  &  & \\
& D_{2\ast} & 0 &  &  &  &  &  & \\
&  & D_{1\ast} & 0 &  &  &  &  & \\
&  &  & D_{\ast} & \left(  T\right)  &  &  &  & \\
\cdots & q_{2} & q_{1} & X & \Delta & 0 &  &  & \\
&  &  &  &  & I_{\mathcal{D}} & 0 &  & \\
&  &  &  &  &  & I_{\mathcal{D}} & 0 & \\
&  &  &  &  &  &  & \ddots & \ddots
\end{array}
\right]  ,
\]
and $\sigma$ is represented as
\[
\sigma=\left[
\begin{array}
[c]{ccccccccc}%
\ddots &  &  &  &  &  &  &  & \\
& \hat{\pi}_{2} &  &  &  &  &  &  & \\
&  & \hat{\pi}_{1} &  &  &  &  &  & \\
&  &  & \hat{\pi} &  &  &  &  & \\
&  &  &  & \left(  \pi\right)  &  &  &  & \\
&  &  &  &  & \rho_{1} &  &  & \\
&  &  &  &  &  & \rho_{1}\circ\alpha &  & \\
&  &  &  &  &  &  & \rho_{1}\circ\alpha^{2} & \\
&  &  &  &  &  &  &  & \ddots
\end{array}
\right]  \text{.}%
\]

\begin{remark}
\label{TransferUnique3}Finally, we note that if $\tau$ is a transfer operator
for $\alpha$, then a contractive covariant representation $(\pi,T)$ has a
unique dilation (up to unitary equivalence) $(\sigma,U)$ where $U$ is unitary
and is adapted to $\tau$ in the sense that $(\sigma,U)$ is the minimal
isometric dilation of the essentially unique adapted coisometric extension
$(\rho, V)$ of $(\pi,T)$.
\end{remark}


\begin{thebibliography}{9}                                                                                                %


\bibitem {wA69}Wm. Arveson, \emph{Subalgebras of $C^{\ast}$-algebras}, Acta
Math. \textbf{123} (1969), 141--224.

\bibitem {vB00}V. Baladi, \emph{Positive transfer operators and decay of
correlations} in Advanced Series in Nonlinear Dynamics, 16. World Scientific
Publishing Co., Inc., River Edge, NJ, 2000. x+314 pp.

\bibitem {rE03}R. Exel, \emph{A new look at the crossed-product of a $C^{*}%
$-algebra by an endomorphism}, Ergodic Theory Dynam. Systems \textbf{23}
(2003), 1733--1750.

\bibitem {MM83}M. McAsey and P. Muhly, \emph{Representations of nonselfadjoint
crossed products}, Proc. London Math. Soc. (3) \textbf{47} (1983), 128--144.

\bibitem {MS98}P. Muhly and B. Solel, \emph{Tensor algebras over $C^{\ast}%
$-correspondences: representations, dilations, and $C^{\ast}$-envelopes}, J.
Funct. Anal. \textbf{158} (1998), 389--457.

\bibitem {MS02}P. Muhly and B. Solel, \emph{Quantum Markov Processes
(Correspondences and Dilations)}, Internat. J. Math. \textbf{13} (2002), 863--906.

\bibitem {pS93}P. Stacey, \emph{Crossed products of $C^{\ast}$-algebras by
$\ast$-endomorphisms}, J. Austral. Math. Soc. Ser. A \textbf{54} (1993), 204--212.

\bibitem {SzNF70}B. Sz.-Nagy and C. Foia\c{s}, \emph{Harmonic Analysis of
Operators in Hilbert Space}, North-Holland, Amsterdam, 1970.


\end{thebibliography}
\end{document}